\documentclass{article}

\usepackage{PRIMEarxiv}

\usepackage[utf8]{inputenc} % allow utf-8 input
\usepackage[T1]{fontenc}    % use 8-bit T1 fonts
\usepackage{hyperref}       % hyperlinks
\usepackage{url}            % simple URL typesetting
\usepackage{booktabs}       % professional-quality tables
\usepackage{amsfonts}       % blackboard math symbols
\usepackage{nicefrac}       % compact symbols for 1/2, etc.
\usepackage{microtype}      % microtypography
\usepackage{lipsum}
\usepackage{fancyhdr}       % header
\usepackage{graphicx}       % graphics
\graphicspath{{media/}}     % organize your images and other figures under media/ folder

\usepackage{tikz}
\usetikzlibrary{arrows,automata,positioning}

\usepackage{mathtools}
\usepackage{xspace}

\usepackage{hyperref}

\usepackage[utf8]{inputenc} % allow utf-8 input
\usepackage[T1]{fontenc}    % use 8-bit T1 fonts
\usepackage{hyperref}       % hyperlinks
\usepackage{url}            % simple URL typesetting
\usepackage{booktabs}       % professional-quality tables
\usepackage{amsfonts}       % blackboard math symbols
\usepackage{nicefrac}       % compact symbols for 1/2, etc.
\usepackage{microtype}      % microtypography
\usepackage{xcolor}         % colors
\usepackage{answers}
\usepackage{amsmath}		% math library

\usepackage[utf8]{inputenc}
\usepackage{amssymb}
\DeclareUnicodeCharacter{211D}{\ensuremath{\mathbb{R}}}
\usepackage[ruled,vlined]{algorithm2e}

\title{Operator-free Equilibrium on the Sphere
}

\author{
  Xiongming Dai  \\
  Division of Computer Science and Engineering \\
  Louisiana State University \\
  Baton Rouge,LA70803, USA\\
  \texttt{\{xdai2\}@email} \\
   \And
  Gerald Baumgartner \\
  Division of Computer Science and Engineering \\
  Louisiana State University \\
  Baton Rouge,LA70803, USA\\
  \texttt{\{gb\}@email} \\
}

\begin{document}
\maketitle

\begin{abstract}
%Particle Gibbs with Ancestor Sampling
We propose a generalized minimum discrepancy, which derives from Legendre's ODE and spherical harmonic theoretics to provide a new criterion of equidistributed pointsets on the sphere. A continuous and derivative kernel in terms of elementary functions is established to simplify the computation of the generalized minimum discrepancy. We consider the deterministic point generated from Pycke's statistics to integrate a Franke function for the sphere and investigate the discrepancies of points systems embedding with different kernels. Quantitive experiments are conducted and the results are analyzed. Our deduced model can explore latent point systems, that have the minimum discrepancy without the involvement of pseudodifferential operators and Beltrami operators, by the use of derivatives. Compared to the random point generated from the Monte Carlo method, only a few points generated by our method are required to approximate the target in arbitrary dimensions.

\end{abstract}

% keywords can be removed

\section{Introduction}
%2003-an efficient algorithm for constructing optimal design of computer experiment
%2006—Book_Design and Modeling for Computer Experiments.... for the second paragraph.
%2010-Generalized Latin Hypercube Design for computer experiment.
%2015-sequential exploration of complex surfaces using minimum energy designs.
%2015-space-filling designs for computer experiments
%2016-GU_dissertation_2016 Minimum Energy designs_extensions
%2017-deterministic sampling of expensive posteriors
%uniquely should be added!
%Step 1: Introduce your topic
%Step 2: Describe the background
%Step 3: Establish your research problem
%Step 4: Specify your objective(s)
%Step 5: Map out your paper
 Quantifying the criterion of equidistributed pointsets on a sphere is of practical importance in numerical analysis \cite{hlawka1982gleichverteilung,grabner1991erdos,grabner1993spherical,rakhmanov1994minimal,rakhmanov1995electrons,kuijlaars1998asymptotics,saff1997distributing,andrievskii1999discrepancy,gotz2001note,bajnok2015constructive,damelin2003energy,sloan2004extremal,hesse2006lower,narcowich2010leveque,brauchart2008optimal,qiu2000numerical},  geophysical \cite{nigrini2007benford,gustafsson2000reconstructing,druskin2010adaptive}, geodetic sciences \cite{cui1997equidistribution,freeden1999constructive,freeden2017integration} and statistics \cite{watson1967another,beran1968testing,gine1975invariant,pycke2007decomposition,pycke2007u,prentice1978invariant}. The advantage of equidistributed point systems is that they are well separated and sufficiently covered such that only a few points are required to approximate the integral. The uniqueness of these points, compared to random points generated from the Monte Carlo method, makes them extensively used in downsampling methods for machine learning. 

For earlier researchers, Freeden obtained explicit identities for the error terms in cubature formulas from the deduction of Green's functions with respect to the Laplace Beltrami operator on the sphere \cite{freeden1980integral}.  {\color{black}Cui} and Freeden extended it further and proposed a generalized discrepancy associated with pseudodifferential operators in $\mathbb{R}^3$ \cite{cui1997equidistribution}. This approach is limited in that the generated point system is only with the kernel-self and cannot explore latent point systems derived from its derivatives further, with mild assumptions.
 
 The purpose of this paper is to study a set of formulas that combines the advantage of Legendre's ODE and further explore latent point systems within error bounds. We consider the properties of the kernel with continuity and derivative, Legendre's ODE and spherical harmonic theoretics to find a new criterion of equidistributed pointsets where the discrepancy becomes smaller, and propose a generalized minimum discrepancy. Our kernel derivative model can explore latent potential point systems that have the minimum discrepancy with operators-free. Our auxiliary intermediaries are spherical harmonic approaches and potential theoretics. 
 
 The paper is divided into three parts. In Section 2, we first introduce a brief abstract of spherical harmonics \cite{muller2006spherical} and the kernel representation for pseudodifferential operators in $\mathbb{R}^3$ \cite{ruzhansky2009pseudo}. For the error estimation of the pointsets, we obtain the upper bound with different order of derivatives of Legendre polynomial and further develop the concept of generalized minimum discrepancy in Section 3.
 
 Our investigation exhibits that, to obtain small discrepancies, point systems on the sphere can be generated from the use of derivatives of kernels without the involvement of the pseudodifferential operators. For different kernels, if they are differentially associated, we can create a mapping $f\in \mathbb{L}^2(\mathbb{S}^d)$ for the pseudodifferential operators. Thus, the generalized minimum discrepancy can be used to reversely deduct the associated pseudodifferential operators.

For certain pseudodifferential operators, we find a closed-form expression of elementary function and group them into different families. It is shown that the two measures of design quality from the point system generated by the generalized minimum discrepancy and the one by minimizing the energy are equivalent. We use the kernel to develop the relation between the points system generated by the minimum energy model and the generalized discrepancies in Section 4. In Section 5, we use the associated kernel to integrate a Franke function for the sphere such that the minimum discrepancies can be obtained under different orders of derivatives, we statistically analyze the discrepancies for different numbers of nodes, and the smoothing parameter estimation for different kernels. we further conduct the experiments of point systems for different kernels on the sphere and compute the discrepancy from the minimum energy perspective. All the tests where the discrepancy of the pointset generated by our methods becomes smaller are valid. The summary of our contributions is outlined in Section 6.
\section{Prerequisites}
%\section{Theory of Spherical Harmonics}
%参考：Discrepancies of Point Sequences on the Sphere and Numerical Integration
{\textbf{Theory of spherical harmonics.}} We use $(x,y,z)$ to represent the element of the three-dimensional Euclidean space $\mathbb{R}^3$ and the Greek alphabet $\xi$ and $\eta$ to represent the vectors of the unit sphere $\mathbb{S}^d$ in $\mathbb{R}^3$. \textbf{x}=$\left\{x_1,...,x_N\right\}$ represents the point system. $\Delta^{\ast}$ represents the Beltrami operator on the unit sphere. A function $f:\mathbb{S}^d \mapsto \mathbb{R}$ possessing $k$ continuous derivatives on $\mathbb{S}^d$ is said to be of the class $C^k(\mathbb{S}^d)$. $C(\mathbb{S}^d)=C^0(\mathbb{S}^d)$ is the class of real continuous scalar-valued functions on $\mathbb{S}^d$. By $\mathbb{L}_2(\mathbb{S}^d)$ we denote the space of Lebesgue square-integrable scalar functions on $\mathbb{S}^d$. Let $Y_{i,j}:i=0,...,n;j=1,...,Z(d,n)$ to be an orthonormalized basis of $\mathbb{L}_2(\mathbb{S}^d)$, where $i$ is called degree, $j$ is the order of the spherical harmonics. The dimension of the space $V_{i}$ of spherical harmonics of order $d+1$ on $\mathbb{S}^{d}$ will be denoted by
\begin{equation}{\label{a1}}
Z(d,i)=(2i+d-1)\frac{\Gamma(i+d-1)}{\Gamma(d)\Gamma(i+1)},
\textbf{1}_{n \gg d}\cdot Z(d,n)= \frac{2}{\Gamma(d)}n^{d-1}.
\end{equation}
The space $V_{i}$ is considered as the eigenspace of the Laplace-Beltrami operator on $\mathbb{S}^{d}$ for the eigenvalue $\lambda_i=-i(i+d-1)$.

The well-known Legendre addition theorem states \cite{muller2006spherical}
\begin{equation}{\label{a2}}
 \sum_{j=1}^{Z(d,i)}Y_{i,j}(\xi)Y_{i,j}(\eta)=\frac{Z(d,i)}{c_{d}}P_i(\xi\cdot \eta), \ \  \xi,\eta\in \mathbb{S}^d,
\end{equation}
where $P_i(x)$ is the Legendre polynomial, an infinitely differentiable eigenfunction of the Legendre operator, orthogonal on the $x\in[-1,1]$ with respect to $(1-x^2)^{d/2-1}$, and it satisfies $P_n(1)=1$, $P_n(x)\le 1$ and $|P_n^{'}(x)|\le \frac{n(n+1)}{2}$. The constant $c_{d}$ denotes the surface area of $\mathbb{S}^d$. 
%Each spherical harmonic $Y_{i,j}$ of exact degree $i$ is an eigenfunction of the negative Laplace-Beltrami operator $\Delta^{\ast}$ with eigenvalue $\lambda_i=-i(i+d-1)$.

\textbf{Functional and distributional spaces.} We consider the space \cite{cui1997equidistribution}
\begin{equation}{\label{a33}} 
H^s{(\mathbb{S}^d)}=\left\{ {f\in C^{\infty}(\mathbb{S}^d)|\sum_{i=0}^{\infty}\sum_{j=1}^{Z(d,i)} f_{i,j}\cdot \hat{i}^{2s}\textless \infty }   \right\},
\end{equation}
where
\[\hat{i}=\left\{\begin{matrix}
1, &\text{if \ \ $i=0$; }\\
i,&\text{otherwise.}
\end{matrix}\right.\]
Then the union of the normalized $Y_{i,j}$ for all $i \in \mathbb{R}$ forms a complete orthonormal system in $\mathbb{L}^2{(\mathbb{S}^d)}$. Thus for $f\in \mathbb{L}^2(\mathbb{S}^d)$, it can be formulated as a Fourier series
\begin{equation}{\label{a3}} f=\sum_{i=0}^{\infty}\sum_{j=1}^{Z(d,i)}\hat{f}_{i,j}Y_{i,j}(\xi),
\end{equation}
where the Fourier coefficients $\hat{f}_{i,j}$ are given by
\begin{equation}{\label{a4}} 
\hat{f}_{i,j}=(f,Y_{i,j})_{\mathbb{L}_2(\mathbb{S}^d)}=\int_{\mathbb{S}^d}^{}f(\xi)Y_{i,j}(\xi)d\sigma_d(\xi),
\end{equation}
satisfying
\begin{equation}{\label{a422}} 
\sum_{i=0}^{\infty}\sum_{j=1}^{Z(d,i)}(1-\lambda_i)^s\left|\hat{f}_{i,j}  \right|\textless \infty,
\end{equation}
where $\sigma_d(\xi)$ denotes the normalized Hausdorff surface measure on the unit sphere $\mathbb{S}^d$ in ${\mathbb{R}}^{d+1}$.

The corresponding inner product in the $H^s(\mathbb{S}^d)$ is
\begin{equation}{\label{a422s}} 
\left\langle f,g \right\rangle_{H^s(\mathbb{S}^d)}=\sum_{i=0}^{\infty}\sum_{j=1}^{Z(d,i)}f_{i,j}g_{i,j}\hat{i}^{2s}, \text{and} \ \  \left\| f \right\|_{H^s(\mathbb{S}^d)}=\sqrt{\sum_{i=0}^{\infty}\sum_{j=1}^{Z(d,i)}f_{i,j}^2\hat{i}^{2s}} \textless \infty.
\end{equation}
From the Cauchy-Schwarz inequality, we obtain
\[
\left( \sum_{i=0}^{\infty}\sum_{j=1}^{Z(d,i)}\left| \hat{f}_{i,j}Y_{i,j}(\xi) \right| \right)^2\le 
\sum_{i=0}^{\infty}\sum_{j=1}^{Z(d,i)}\left| \hat{f}_{i,j}^2\hat{i}^{2s} \right|\cdot
\sum_{i=0}^{\infty}\sum_{j=1}^{Z(d,i)}\left| Y_{i,j}^2\hat{i}^{2s} \right|= \left\| f \right\|_{H^s(\mathbb{S}^d)}\sum_{i=0}^{\infty}Z(d,i)i^{-2s}.
\]
As $Z(d,i) \le i^{d-1}$, thus the series uniformly converges for $d-1-2s<-1  \Rightarrow  s > \frac{d}{2}$.

Thus, the spherical harmonic expansion of any function $f$ in $H^s(\mathbb{S}^d)$ will converge uniformly for $s > \frac{d}{2}$. This is significant since there are functions in $C^k(\mathbb{S}^d)$ which do not allow a uniformly convergent for spherical harmonic series \cite{cui1997equidistribution,gronwall1914degree}. For our experiment in Section 5, we use $s>2$.

\textbf{Pseudodifferential operator.} $H^s(\mathbb{S}^d)\subset C^k(\mathbb{S}^d),$ for $s > \frac{d}{2}$. Let $\left\{ A_i \right\}_{i\in\mathbb{R}^+}$ be a sequence of real numbers $A_i$ satisfying
\[
\lim_{i \to 0} \frac{\left| A_i \right|}{(i+\frac{d-1}{2})^{\alpha}}=\text{const}\neq 0
\]
for a certain $\alpha \in \mathbb{R}^+$. Then a pseudodifferential operator of order $\alpha$, \textbf{A} from $H^s(\mathbb{S}^d)$ to $H^s(\mathbb{S}^{d-\alpha})$ is defined by
\begin{equation}\label{452 3d}
\textbf{A}f=\sum_{i=0}^{\infty}\sum_{j=1}^{Z(d,i)}A_i\hat{f}_{i,j}Y_{i,j}(\xi),f\in H^s(\mathbb{S}^d). 
\end{equation}
The sequence $\left\{ A_i \right\}_{i\in\mathbb{R}^+}$ is called the spherical symbol of $\textbf{A}$. It is obvious that, for a pseudodifferential operator $\textbf{A}$ of order $s$, equation ~\eqref{a33} $H^s(\mathbb{S}^d)$ can be equivalently expressed as 
\[
H^s(\mathbb{S}^d)=\left\{ f:\mathbb{S}^d \to \mathbb{R} |\textbf{A}f\in \mathbb{L}_2(\mathbb{S}^d)\right\}.
\]
The relation between the pseudodifferential operator $\textbf{A}$ on the sphere and the Beltrami operator $\Delta^{\ast}$ for a certain elementary functional representation is provided by \cite{cui1997equidistribution}. We consider equation ~\eqref{a422s}, the kernel $K$ associated with the space $H^s(\mathbb{S}^d)$ and the inner product $\left\langle f,g \right\rangle_{H^s(\mathbb{S}^d)}$ is
\begin{equation}\label{a3ad}
K(\xi\cdot \eta) =\sum_{i=0}^{\infty}\sum_{j=1}^{Z(d,i)} \frac{1}{\hat{i}^{2s}} \cdot Y_{i,j}(\xi) \cdot Y_{i,j}(\eta)=\sum_{i=0}^{\infty}\frac{Z(d,i)}{\hat{i}^{2s} \cdot c_d} \cdot P_i(\xi \cdot \eta),
\end{equation}
for invariant pseudodifferential operator \textbf{A} on the sphere, it can be simplified into 
\begin{equation}\label{a3a}
K_\textbf{A}(\xi\cdot \eta) =\sum_{i=0}^{\infty}\sum_{j=1}^{2n+1} A_n \cdot Y_{i,j}(\xi) \cdot Y_{i,j}(\eta)=\sum_{i=0}^{\infty}\frac{2n+1}{4\pi} \cdot A_n \cdot P_n(\xi \cdot \eta).   
\end{equation}
The equation ~\eqref{a3a} can be further simplified by convolution into 
\[\textbf{A}f=K_{\textbf{A}}\ast f=\int_{H^s(\mathbb{S}^d)}K_\textbf{A}(\xi\cdot \eta)f(\xi)d\sigma_d(\xi).
\]
The kernel $K_\textbf{A}(\xi\cdot \eta)\in H^{-(\alpha+\zeta)}(\mathbb{S}^d) $  for all $\zeta>0$ \cite{cui1997equidistribution}.
\section{Operator-free Equilibrium by Derivatives}
In this section, we focus on the discrepancies of equilibrium from different self-joint kernels.
%Let $x_i \in \mathbb{R}^3,(i=1,...,N)$ be Cartesian coordinates of a Euclidean space.
The problem can be stated as follows: There exist coefficients $a_i$ such that $\sum_{i=1}^{N}a_if(x_i)$ is a good approximation to $\frac{1}{4\pi}\int_{\mathbb{S}^d}f(x)d\omega(x)$ in a certain upper bound for any $f\in \mathbb{L}^2(\mathbb{S}^d)$, as $N\to \infty$.

\textbf{Theorem 1}  Let $\textbf{A}$ be a pseudodifferential operator of order $s$, $s>1$, with the symbol ${A_n}$ satisfying  $A_n>0,n \geq 1$. Let $m$ denote the order of the highest derivative
for Legendre polynomial $P_n(t)$, for any function $\textbf{A}f(x)\in \mathbb{L}^2(\mathbb{S}^d)$ and $m \le N, m\in\mathbb{R}$, we have the estimate
\begin{equation}\label{a345}
\left| \sum_{i=1}^{N}a_if(x_i)-\frac{1}{4\pi}\int_{\mathbb{S}^d}^{}f(x)d\omega(x) \right|\le \frac{1}{N}\sqrt{\left[ \sum_{t=1}^{N}\sum_{i=1}^{N}\sum_{n=1}^{\infty} \frac{Z(d,i)}{A_n^2} \frac{\partial^mP_n}{(\partial (\eta_i \cdot \eta_t))^m}(\eta_i \cdot \eta_t)\right]}\left\| \textbf{A}f(x) \right\|_{L^2}.  
\end{equation}
\textbf{Proof} From ~\eqref{a33}, we can induce $f(\xi)\in C^{\infty}(\mathbb{S}^d)$. As $s>1$ and $d \geq 2$, the spherical harmonic expansion of any function $f(\xi)\in H^s(\mathbb{S}^d) $ will converge uniformly, we have
\begin{equation}\label{a3451}
f(\xi) =\sum_{n=0}^{\infty}\sum_{j=1}^{2n+1}f_{n,j}Y_{n,j}(\xi),\xi\in \mathbb{S}^d.  
\end{equation}
We discrete the surface with $d\omega(\eta)$ on the sphere. From \cite{freeden1980integral}, we get
\begin{equation}\label{a3452}
f(\xi) =\frac{1}{4\pi}\int_{\mathbb{S}^d}^{}f(\eta)d\omega(\eta)-\int_{\mathbb{S}^d}^{}\sum_{k=0}^{\infty}\frac{1}{k(k+1)-\lambda_k}\sum_{j=1}^{2k+1}Y_{k,j}(\xi)Y_{k,j}(\eta)\Delta_{\xi}^*f(\eta)d\omega(\eta).  
\end{equation}
Given $\xi=\eta_i$, $i\in [1,N]$ and $Y_{k,j}(\xi)= \frac{1}{N}\sum_{i=1}^{N}Y_{k,j}(\eta_i)$, it leads to
\begin{equation}\label{a3453}
\sum_{i=1}^{N}a_if(\eta_i) =\frac{1}{4\pi}\int_{\mathbb{S}^d}^{}f(\eta)d\omega(\eta)-\frac{1}{N}\int_{\mathbb{S}^d}^{}\sum_{k=0}^{\infty}\frac{1}{k(k+1)-\lambda_k}\sum_{j=1}^{2k+1}\sum_{i=1}^{N}Y_{k,j}(\eta_i)Y_{k,j}(\eta)\Delta_{\xi}^*f(\eta)d\omega(\eta). 
\end{equation}
From the Cauchy-Schwarz inequality and the Legendre addition theorem, we get \cite{cui1997equidistribution}
\begin{equation}
\label{nined1s}
\begin{split}
&\left| \sum_{i=1}^{N}a_if(\eta_i)-\frac{1}{4\pi}\int_{\mathbb{S}^d}^{}f(\eta)d\omega(\eta) \right|\\
&\le \frac{1}{N}\sum_{k=0}^{\infty}\sum_{j=1}^{2k+1}\sum_{i=1}^{N}\frac{1}{k(k+1)-\lambda_k}\int_{\mathbb{S}^d}^{}Y_{k,j}(\eta_i)Y_{k,j}(\eta)\Delta_{\xi}^*f(\eta)d\omega(\eta)\\
&=\frac{1}{N}\int_{\mathbb{S}^d}^{}f(\eta)\sum_{k=0}^{\infty}\sum_{j=1}^{2k+1}\sum_{i=1}^{N}\frac{\Delta_{\xi}^*}{k(k+1)-\lambda_k}Y_{k,j}(\eta_i)Y_{k,j}(\eta)d\omega(\eta)\\
&= \frac{1}{N}\sqrt{\int_{\mathbb{S}^d}^{}f^2(\eta)d\omega(\eta)}\cdot \sqrt{\int_{\mathbb{S}^d}^{}\left(\sum_{k=0}^{\infty}\sum_{j=1}^{2k+1}\sum_{i=1}^{N}\frac{Y_{k,j}(\eta_i)Y_{k,j}(\eta)}{A_n}\right) ^2d\omega(\eta)}\\
&=\frac{1}{N}\left\| \textbf{A}f(\xi) \right\|_{L^2}\sqrt{\sum_{n=1}^{\infty} \sum_{j=1}^{2n+1}\left( \frac{\sum_{i=1}^{N}Y_{n,j}(\eta_i)}{A_n} \right)^2}\\
&=\frac{1}{N}\left\| \textbf{A}f(\xi) \right\|_{L^2}\sqrt{ \sum_{n=1}^{\infty} \sum_{j=1}^{2n+1}\sum_{t=1}^{N}\sum_{i=1}^{N} \frac{Y_{n,j}(\eta_i)Y_{n,j}(\eta_t)}{A^2_n} }\\
&=\frac{1}{N}\left\| \textbf{A}f(\xi) \right\|_{L^2}\sqrt{\sum_{n=1}^{\infty} \sum_{t=1}^{N}\sum_{i=1}^{N} \frac{2n+1}{4\pi A^2_n}P_n(\eta_i \cdot \eta_t) }.
\end{split}
\end{equation}
Here, we first focus on the deduction of Legendre polynomial $P_n(\eta_i\cdot \eta_t)$ recurrence relations. Differentiating the generating function \cite{arfken1972mathematical}
\begin{equation}\label{a321}
g(x,t)=(1-2xt+t^2)^{-\frac{1}{2}}=\sum_{n=0}^{\infty}P_n(x)t^n,\left| t \right|<1, 
\end{equation}
with respect to $x$, we get
\begin{equation}\label{a321a}
\frac{\partial g(x,t)}{\partial x}=\frac{t}{(1-2xt+t^2)^\frac{3}{2}}=\sum_{n=0}^{\infty}P^{'}_n(x)t^n.
\end{equation}

Substituting  ~\eqref{a321} to  ~\eqref{a321a}, we get
\begin{equation}\label{a321a1}
(1-2xt+t^2)\sum_{n=0}^{\infty}P^{'}_n(x)t^n-t\sum_{n=0}^{\infty}P_n(x)t^n=0,
\end{equation}
which leads to
\begin{equation}\label{a321a1a}
P^{'}_{n+1}(x)+P^{'}_{n-1}(x)=2xP^{'}_{n}(x)+P_n(x).
\end{equation}
Differentiating the following Bonnet's recursion formula
\begin{equation}\label{a321a1a4}
(2n+1)xP_n(x)=(n+1)P_{n+1}(x)+nP_{n-1}(x),
\end{equation}
with respect to $x$, and adding 2 times $\frac{d}{dx}$~\eqref{a321a1a4} to $(2n+1)$ times ~\eqref{a321a1a}, we get
\begin{equation}\label{aa4}
(2n+1)P_n(x)=P^{'}_{n+1}(x)-P^{'}_{n-1}(x).
\end{equation}
From the above, we can also find that
\begin{equation}\label{ab4}
P^{'}_{n+1}(x)=(2n+1)P_n(x)+(2(n-2)+1)P_{n-2}(x)+(2(n-4)+1)P_{n-4}(x)+\cdots, \end{equation}
or equivalently
\begin{equation}\label{ab24}
P^{'}_{n+1}(x)=\frac{2}{\left\| P_n \right\|^2}P_n(x)+\frac{2}{\left\| P_{n-2} \right\|^2}P_{n-2}(x)+\frac{2}{\left\| P_{n-4} \right\|^2}P_{n-4}(x)+\cdots,
\end{equation}
where $\left\| P_n(x) \right\|$ is the norm over the interval $x \in[-1,1]$
\begin{equation}\label{ab241}
\left\| P_n \right\|=\sqrt{\int_{-1}^{1}(P_n(x))^2dx}=\sqrt{\frac{2}{2n+1}},
 \end{equation}
satisfying from Rodigue's formula
\begin{equation}\label{ab241a}
P_n(x)=\frac{1}{2^nn!}\frac{d^n}{dx^n}(x^2-1)^n.
\end{equation}
The standardization $P_n(1)=1$ fixes the normalization of the Legendre polynomials, since they are also orthogonal with respect to the same norm, and can be recursively nested to the order of the highest derivative $m$ from Equations ~\eqref{ab24} and ~\eqref{ab241a}, we can find that there exists a $m$, $m \le N,m \in \mathbb{R}$ satisfying $P^{(m)}_{n}(x)=\sum_{n=1}^{N}\beta_nP_n(x)$ for a certain $\beta_n\in \mathbb{R}$, thus ~\eqref{nined1s} can be rewritten as
\begin{equation}
\label{nined1s5}
\begin{split}
&\left| \sum_{i=1}^{N}a_if(\eta_i)-\frac{1}{4\pi}\int_{\mathbb{S}^d}^{}f(\eta)d\omega(\eta) \right|\\
&\le \frac{1}{N}\left\| \textbf{A}f(\xi) \right\|_{L^2}\sqrt{ \sum_{n=1}^{\infty} \sum_{j=1}^{2n+1}\sum_{i=1}^{N} \frac{2n+1}{4\pi A^2_n}P_n(\eta_i \cdot \eta_t) }\\
&=\frac{1}{N}\left\| \textbf{A}f(\xi) \right\|_{L^2}\sqrt{ \sum_{n=1}^{\infty} \sum_{t=1}^{N}\sum_{i=1}^{N} \frac{2n+1}{4\pi A^2_n}\frac{\partial^mP_n }{(\partial (\eta_i \cdot \eta_t))^m}(\eta_i \cdot \eta_t) }.
\end{split}
\end{equation}
This completes the proof.

 Theorem 1 shows that the error highly depends on the pointset. This gives rise to the definition of generalized minimum discrepancy.
 
 \textbf{Generalized minimum discrepancy.} Let $\textbf{A}$ be a pseudodifferential operator of order $s$, $s>1$, with symbol ${A_n}$, $A_n\neq 0$ for $n \geq 1$. Then the generalized minimum discrepancy associated with a pseudodifferential operator $\textbf{A}$ is defined by
\begin{equation}\label{a3456}
D_{\min}(\textbf{x};\textbf{A})=\min(\frac{1}{N}\sqrt{\sum_{t=1}^{N}\sum_{i=1}^{N}\sum_{n=1}^{\infty} \frac{Z(d,i)}{A_n^2} \frac{\partial^mP_n }{(\partial (\eta_i \cdot \eta_t))^m}(\eta_i \cdot \eta_t)}),   
\end{equation}
where $m\in[0,N]$ denotes the order of the highest derivative
for Legendre polynomial $P_n(\cdot)$.

The minimum discrepancy shows that, for $m\in [0,N]$, there exist, different groups of point sets, where the minimum discrepancy of the asymptotically distributed is the $m$-th order of the highest derivatives. This can be interpreted intuitively as follows: Given a point set $\textbf{x}$ on the sphere $\mathbb{S}^d$, the measure for the quality of the distribution is the spherical cap discrepancy
\begin{equation}\label{a33f}
D(\textbf{x})=\underset{C\subseteq \mathbb{S}^d}{\text{sup}}\left| \frac{1}{N}\sum_{i=1}^{N} \delta_C(x_i)-\frac{1}{4\pi} f_C(\xi)d\omega(\xi)\right|,  
\end{equation}
where the supremum ranges over all spherical caps $C\subseteq \mathbb{S}^d$ (intersections of ball and $\mathbb{S}^d$) and $\delta_C$ represent the Dirac delta measure that associates to $C$.
The discrepancy simply measures the maximal deviation between the discrete distribution $\textbf{x}$ and the normalized surface measure. Let $f(\xi) \in H^s(\mathbb{S}^d),s>1$, we have 
\begin{equation}
\label{nined1s53}
\begin{split}
&D(\textbf{x})=\underset{C\subseteq \mathbb{S}^d}{\text{sup}}\left| \frac{1}{N}\sum_{i=1}^{N} \delta_C(x_k)-\frac{1}{4\pi} f_C(\xi)d\omega(\xi)\right|\\
&\approx \left| \frac{1}{N}\sum_{i=1}^{N}\delta_C(x_i)-\frac{1}{4\pi}\int_{\mathbb{S}^d}f(\xi)d\omega(\xi)^{} \right|\\
&=\frac{1}{N}\left\| \textbf{A}f(\xi) \right\|_{L^2}\sqrt{\sum_{n=1}^{\infty} \sum_{t=1}^{N}\sum_{i=1}^{N} \frac{2n+1}{4\pi A^2_n}\frac{\partial^mP_n }{(\partial (\eta_i \cdot \eta_t))^m}(\eta_i \cdot \eta_t)}\\
&\le \frac{1}{N}\left\| \textbf{A}f(\xi) \right\|_{L^2}\sqrt{\sum_{n=1}^{\infty} \sum_{t=1}^{N}\sum_{i=1}^{N} \frac{2n+1}{4\pi A^2_n}P^{}_{n}(\eta_i \cdot \eta_t)}.
\end{split}
\end{equation}
Comparing to \cite{cui1997equidistribution}, we consider the different directions of the point $x_i$ on the sphere by the derivatives with $m$ order($P^{(m)}_{n}(\eta_i \cdot \eta_t)$), not limited only one direction
where $m=0$ with $P^{}_{n}(\eta_i \cdot \eta_t)$. Thus, the generalized minimum discrepancy exhibits a more wide range of exploring the candidates to asymptotically distribute the spherical cap $C$.

\textbf{Lemma 1} Let \textbf{A}, \textbf{B} be two pseudodifferential operators of order $s_1,s_2(s_1>1,s_2>1),$ and with symbols $\left\{A_n\right\}$,$\left\{B_n \right\}$ satisfying $A_n> 0,B_n> 0$ for $n \ge 1$, respectively. $K_A(\xi \cdot \eta)$ and $K_B(\xi \cdot \eta)$ satisfying ~\eqref{a3a}. If
\begin{equation}\label{a33f1}
(-1)^nc_n\frac{\partial^nK_A }{(\partial (\xi \cdot \eta))^n}(\xi \cdot \eta)=K_B(\xi \cdot \eta),n\in \mathbb{R}^+, 
\end{equation}
with the factor $c_0=1, c_n=\frac{1}{(n-1)!}$, there exists a $f\in \mathbb{L}^2(\mathbb{S}^d)$, such that $B_n=f(A_n)$ and $D_{\min}(\textbf{x};\textbf{A})=D_{\min}(\textbf{x};\textbf{B})$. We call the discrepancies from both $\left\{A_n\right\}$, and $\left\{B_n\right\}$ belong to the same family discrepancies, the associated kernels $K_A$ and $K_B$ belong to the same family kernel.

\textbf{Proof} From ~\eqref{a3a}, $K_A(\eta \cdot \xi)\propto \sum_{n=0}^{\infty}P_n(\eta \cdot \xi)$, as $P_n(\xi \cdot \eta)$ is normalized orthogonal basis, from the Rodigue's formula ~\eqref{ab241a}, it is the $n$-th order derivative in $ \left[ -1,1 \right]$, we obtain  \begin{equation}\label{a3ag}
\frac{\partial^{m_a}K_A }{(\partial (\xi \cdot \eta))^{m_a}}(\xi \cdot \eta) =\sum_{n=0}^{\infty}Z(d,i) \cdot A_n \cdot \frac{\partial^{m_a}P_n }{(\partial (\eta_i \cdot \eta_t))^{m_a}}(\eta_i \cdot \eta_t).    
\end{equation}
Substitute ~\eqref{a33f1} to ~\eqref{a3ag}, it yields
\begin{equation}\label{a3agf}
(-1)^nc_n\frac{\partial^{n+m_b}K_A }{(\partial (\xi \cdot \eta))^{n+m_b}}(\xi \cdot \eta) =\sum_{n=0}^{\infty}Z(d,i) \cdot B_n \cdot \frac{\partial^{m_b}P_n }{(\partial (\eta_i \cdot \eta_t))^{m_b}}(\eta_i \cdot \eta_t).    
\end{equation}
From ~\eqref{ab24}, each derivative item on the right-side $\frac{\partial^{m_b}P_n }{(\partial (\eta_i \cdot \eta_t))^{m_b}}(\eta_i \cdot \eta_t)$ can be represented by the normalized basis $P_i(\eta_i \cdot \eta_t)$, with orthogonality and completeness, there exists a piecewise continue function $f(\cdot)\in \mathbb{L}^2(\mathbb{S}^d)$ with finitely many discontinuities in $ \left[ -1,1 \right]$, the sequence of sums 
\begin{equation}\label{ad2456}
f_n(x,A_n)=\sum_{i=0}^{n}a_i\cdot B_i \cdot P_i(x),
\end{equation}
converges in the mean to $f(x,\textbf{A})$ as $ n\rightarrow \infty$, provided we take
\begin{equation}\label{ad23456}
a_i=\frac{2i+1}{2}\int_{-1}^{1}f(x,\textbf{A})P_i(x)dx.
\end{equation}
For pseudodifferential operator $\textbf{A}$, we obtain
\begin{equation}\label{ad456}
D_{\min}(\textbf{x};\textbf{A})=\min(\frac{1}{N}\left[ \sum_{t=1}^{N}\sum_{i=1}^{N}\sum_{n=1}^{\infty} \frac{Z(d,i)}{A_n^2} \frac{\partial^{m_a}P_n }{(\partial (\eta_i \cdot \eta_t))^{m_a}}(\eta_i \cdot \eta_t)\right]^\frac{1}{2}), {m_a}\in[0,N],  
\end{equation}
comparing to $\textbf{B}$
\begin{equation}\label{a34g56}
D_{\min}(\textbf{x};\textbf{B})=\min(\frac{1}{N}\left[ \sum_{t=1}^{N}\sum_{i=1}^{N}\sum_{n=1}^{\infty} \frac{Z(d,i)}{B_n^2} \frac{\partial^{m_b}P_n }{(\partial (\eta_i \cdot \eta_t))^{m_b}}(\eta_i \cdot \eta_t)\right]^\frac{1}{2}), m_b\in[0,N]. \end{equation}
Combine with ~\eqref{a3agf}, it is obvious that
$\frac{\partial^{m_b}P_n }{(\partial (\eta_i \cdot \eta_t))^{m_b}}(\eta_i \cdot \eta_t)\propto \frac{\partial^{m_b+n}P_n }{(\partial (\eta_i \cdot \eta_t))^{m_b+n}}(\eta_i \cdot \eta_t)$. Thus, $D_{\min}(\textbf{x};\textbf{A})=D_{\min}(\textbf{x};\textbf{B})$. This completes the proof.

From lemma 1 we prove that the generalized minimum discrepancy can be used to reversely deduct the associated pseudodifferential operators and find that for different kernels if they are differentially associated, we can create a mapping $f\in \mathbb{L}^2(\mathbb{S}^d)$ for the pseudodifferential operators. Using this property, we can extend the potential theoretics for the logarithmic energy kernel and Riesz kernel. 

\textbf{Equidistribution in $H^s(\mathbb{S}^d)$.} A point system $\textbf{x}$ is called $\textbf{A}$-equidistributed in $H^s(\mathbb{S}^d)$, $s>1$ if the generalized discrepancy associated with a pseudodifferential operator $\textbf{A}$ of order $s$, $s>1$ satisfies
\begin{equation}\label{a3457}
\lim_{N \to \infty} D_{\min}(\textbf{x};\textbf{A})  =0.
\end{equation}
If $\textbf{x}$ is well equidistributed in $H^s{(\mathbb{S}^d)}$, $s>1$. For $s'>s$, we generally need more points such that the point system also uniformly equidistributes in $H^{s'}{(\mathbb{S}^d)}$. Thus, we try to use $s$ as small as possible \cite{cui1997equidistribution}. However, for the computation of ~\eqref{a3456}. We need to calculate the series expansion in terms of Legendre polynomials derivative with order $m$. From ~\eqref{ab241a}, the complexity is $\mathcal{O}(2^n)$. 
It is not applicable to use ~\eqref{a3456} directly for the solver of the generalized minimum discrepancy.

For certain pseudodifferential operators, we can find a closed-form expression for ~\eqref{a3456}, which has been verified by statistics. Combining with ~\eqref{a3a}, we get
\begin{equation}\label{ad56}
D_{\min}(\textbf{x};\textbf{A})\propto \min(\frac{1}{N}\left[ \sum_{i=1}^{N}\sum_{n=1}^{\infty}  \frac{\partial^mK_A }{(\partial (\eta_i \cdot \eta_t))^m}(\eta_i \cdot \eta_t)\right]^\frac{1}{2}), m\in[0,N].  
\end{equation}
Certain kernels with the corresponding complicated statistics are provided in \cite{choirat2013computational}, we provide some complicated cases as follows.

(1) Gin\'e's statistic:
 $K_A(\eta_i,\eta_j)=\frac{1}{2}-\frac{2}{\pi}\sin \cos^{-1}(\eta_i \cdot \eta_j)$. Where $A_n^2=+\infty$ for $n$ odd and $A_n^2=\frac{n-1}{n+2}\cdot(\frac{\Gamma(\frac{n}{2})}{\Gamma(\frac{n+1}{2})})^2$ for $n$ even \cite{gine1975invariant,gradshteyn2014table}. 

(2) Beran's form of Ajne's statistic:
$K_A(\eta_i,\eta_j)=\frac{1}{4}-\frac{1}{2 \pi}\cos^{-1}(\eta_i \cdot \eta_j)$. Where $A_n^2=+\infty$ for $n$ even and $A_n^2=n^2 \cdot (\frac{\Gamma(\frac{n+3}{2})}{\Gamma(\frac{n+2}{2})})^2$ for $n$ odd \cite{gine1975invariant,prentice1978invariant}.

(3) Pycke's statistic:
$K_A(\eta_i,\eta_j)=-\frac{1}{4 \pi}\ln\frac{e}{2}(1-\eta_i \cdot \eta_j)$. Where $A_n^2=n(n+1)$ \cite{pycke2007decomposition,pycke2007u}.

(4) Cui-Freeden Discrepancy:
$K_A(\eta_i,\eta_j)=1-2\ln(1+\sqrt{\frac{1-\eta_i \cdot \eta_j}{2}})$. Where $A_n^2=n(n+1)(2n+1)$ \cite{cui1997equidistribution}.
 
(5) Riesz kernels \cite{damelin2008walk}:
 For $x_i,x_j\in C$ we define

\begin{equation}\label{3wdd}
K(x_1,x_2)=
\begin{cases}
 &\text{sign}(s) \cdot \left\| x_1-x_2 \right\|_2^{-s},s\neq 0, \\ 
 &-\text{ln}\left\| x_1-x_2 \right\|_2^{-2},s= 0, 
 \end{cases}
\end{equation}
 where $\left\| \cdot \right\|_2$ is the Euclidean distance. The logarithmic potential is at the case $s=0$ and the Coulombic potential is at the case $s=1$, respectively. For a unit sphere, we transform it into the vector format as follows \cite{choirat2013computational}:
\begin{equation}\label{3rwdd}
K_A(\eta_i\cdot \eta_j)=
\begin{cases}
 &\text{sign}(s) \cdot \left| 2(1-\eta_i \cdot \eta_j) \right|^{-\frac{s}{2}},s\neq 0, \\ 
 &-\ln2(1-\eta_i \cdot \eta_j),s= 0, 
 \end{cases}
\end{equation}
for $\eta_i \cdot \eta_j \in [-1,1)$. When $s\neq 0$, $A_n^2=\frac{2^{s-2}\Gamma(\frac{s}{2})\Gamma(-\frac{s}{2}+n+2)}{\pi\Gamma(\frac{s}{2}+n)\Gamma(1-\frac{s}{2})}$.
For $s<2$,
$K_A(\eta_i \cdot \eta_j)=\left| 2(1-\eta_i \cdot \eta_j) \right|^{-\frac{s}{2}}-\frac{2^{-s}}{1-\frac{s}{2}}$. When $s=0$, $A_n^2=\frac{n(n+1)}{4\pi}$,
$K_A(\eta_i \cdot \eta_j)=-\ln2(1-\eta_i \cdot \eta_j)-\ln\frac{e}{4}$.
This is a version of Pycke's statistics. Thus, from Lemma 1, the logarithmic potential, Coulombic potential, Pycke's statistics, and Riesz kernel belong to the same family kernel. Thus, we can expand into the generality that for the same family of kernels, we can bypass the pseudodifferential operators from deriving the kernel so as to obtain the minimum discrepancy. We call the minimum discrepancy kernel with more global behavior.

\section{Discrepancy Inequalities via Energy Methods} In physics experiments, we use the principle of mutual repulsion of charges to investigate how to distribute $N$ point charges over a surface $M$, usually by minimizing the sum of all potential energies to obtain the optimal configuration of these charges. The study of the accurate distribution of the charges is the subject of classical potential theory, which shows that the energy integral can be solvable or approximated amongst all Borel probability measures supported on the space. This optimal measure depends highly on the curvature of the position on the surface and the value of $s$ and $d$. 

\textbf{Kernels, energy and measures.}
Let $\Omega$ denote a compact and measurable subset of Euclidean space in $\mathbb{R}^d$ whose $d$-dimensional Borel measure (\textit{charge distributions}) $\mu \subset (\Omega,\mathbb{R}^d)$, is finite, and in the context of energy, $K$ denote a bi-Lipschitz mapping from $\Omega \times \Omega $ to $\mathbb{R}^d$, for a collection of $N(\geq 2)$ distinct points of configuration in $\Omega$, let $X_{1:N}={x_1,...,x_N}$, we define the energy of $X_{1:N}$ to be
\begin{equation}{\label{b2}}
E(X_{1:N}):=\frac{1}{N^2}\sum_{i=1}^{N}\sum_{j=1,j\neq i}^{N}K(x_i,x_j)=\frac{1}{N^2}\sum_{i\neq j}^{}K(x_i,x_j),
\end{equation}
and let 
\begin{equation}{\label{b3}}
\mathcal{E}(\Omega ,N):=\inf\{E(X_{1:N}):X_{1:N}\subset \Omega,\left | X_{1:N} \right |=N  \}
\end{equation}
be the minimal discrete $N$-point energy of the configuration in $\Omega$, where $\left | X_{1:N} \right |$ represents the cardinality of the set $X_{1:N}$.
The measure of the total charge distributed on $\Omega$ can be expressed as $Q(\mu):=\mu(\Omega)=\int_{\Omega}^{}d\mu(x)$. 

For all signed Borel measures (continuous charge distributions) $\mu$ on $\mathbb{S}^d$, the energy integral
\begin{equation}\label{a342}
E(\mu)=\underset{\mathbb{S}^d\times \mathbb{S}^d} {\int\int} K(\xi\cdot \eta)d\mu(\xi)d\mu(\eta) \geq 0, \text{for all} \ \  \mu\neq 0.  
\end{equation}
A measure is a countably additive, non-negative, extended real-valued function
defined on a $\sigma$-algebra $\mathcal{T}$(a nonempty collection of subsets of $X$ closed under complement, countable unions, and countable intersections).

A measure $\mu$ on a measurable space $(X,  \mathcal{T})$ is a mapping
\[ \mu:\mathcal{T} \to [0,\infty]\]
such that  (1) $\mu(\emptyset )=0$; (2) if $\left\{ T_i \in \mathcal{T}:i \in \mathbb{N}\right\}$ is a countable  collection of pairwise disjoint sets in $\mathcal{T}$, then
\[
\mu(\cup_{i=1}^{\infty}T_i)=\sum_{i=1}^{\infty}\mu(T_i).
\]
Let $\delta_x \in(X,\mathcal{T})$ represent the Dirac delta measure that associates to a unit charge at the point $x\in X$, satisfying $\int_{X{'}}d\delta_x(\xi)=1$ for all measurable sets $X{''}\subseteq X$ with $x\in X{''}$. For the empirical distribution of set $X'$, defined as
\begin{equation}\label{a34a}
\mu_{X'}:=\frac{1}{N}\sum_{i=1}^{N}\delta_{x_i},
\end{equation}
we have $\mathcal{E}(X')=\mathcal{E}(\mu_{X'})$.

The quadratic form in ~\eqref{a342} can be used to define the inner product for the charge distribution
\begin{equation}\label{a34as2}
\left\langle \mu,\rho \right\rangle_{(X,\mathcal{T})}=\underset{\mathbb{S}^d\times \mathbb{S}^d} {\int\int} K(\xi\cdot \eta)d\mu(\xi)d\rho(\eta),  
\end{equation}
and the energy then associates with the square norm of the measure
\begin{equation}\label{a3d4a22}
\mathcal{E}(\mu)=\left\| \mu \right\|_{(X,\mathcal{T})}^2. 
\end{equation}
The discrepancy of the measure $\rho$ with respect to the measure $\mu$ is defined as in \cite{damelin2010energy} as
\begin{equation}\label{a3d4a2d2} 
D(\rho;\mu):=\left\| \rho-\mu \right\|_{(X,\mathcal{T})}. 
\end{equation}
Both the energy and the discrepancy highly depend on the choice of the kernel and the charge distribution.

For every signed measure $\mu\in (X,\mathcal{T})$, the potential field induced by the charge distribution by
\begin{equation*}\label{a3d4a2dd2}
f_{\mu}(x)=\int_{\Omega}^{}K(x,y)d\mu(y). 
\end{equation*}
Let $\mathcal{U}(K)$ represent the domain of measures of potential fields, the inner product on $\mathcal{U}(K)$ 
\begin{equation}\label{a34ads2}
\left\langle f_\mu,f_\rho \right\rangle_{\mathcal{U}(K)}=\left\langle \mu,\rho \right\rangle_{(X,\mathcal{T},\mu)}  \ \ \forall f_{\mu},f_{\rho}\in{\mathcal{U}(K)}. 
\end{equation}
The energy can be rewritten into the format with respect to the potential fields by
\begin{equation}\label{a342dd}
E(\mu)=\underset{\mathbb{S}^d\times \mathbb{S}^d} {\int\int} K(\xi\cdot \eta)d\mu(\xi)d\mu(\eta) =\int_{\Omega}^{}f_{\rho}(\xi)d\mu(\xi)=\left\langle f_\mu,f_\rho \right\rangle_{\mathcal{U}(K)}.  
\end{equation}
The energy of the charge distribution can be interpreted physically as the integration of the potential fields against the measure $\mu$.

%\textbf{Relationship between  Minimum Energy  and Generalized Minimum Discrepancy.}
\textbf{Theorem 2} {\color{black}The two} measures of design quality from the point system generated by $D_{\min}(\textbf{x};\textbf{A})$ and the one by minimizing the energy are equivalent, it satisfies $D_{\min}(\textbf{x};\textbf{A})=\sqrt{\mathcal{E}(\mu(\textbf{x})-\omega(\textbf{x}))}$ with $\mu,\omega\in {(X,\mathcal{T})}$.

\textbf{Proof} From ~\eqref{a3456}, if $m=m^*$ is the current order of derivatives where exists the minimum discrepancy, and the corresponding pointset $X^*=\left\{x_1^*,...,x_N^*  \right\}$ on the sphere $\mathbb{S}^d$, the sphere cap discrepancy from ~\eqref{a33f},
\begin{equation}\label{a33eef}
D_{\min}(\left\{X^*; \textbf{A}\right\})=\underset{C\subseteq \mathbb{S}^d}{\text{sup}}\left| \frac{1}{N}\sum_{k=1}^{N} \delta_C(x_k^*)-\frac{1}{4\pi} f_C(\xi)d\omega(\xi)\right|.
\end{equation}
If exists a $\mu(\cdot)$ such that $X^*\sim \mu(\cdot)$, combine 
~\eqref{a342dd}, ~\eqref{a33eef} can be rewritten into
\begin{equation}
\label{nined1ws5w3}
\begin{split}
D_{\min}(\left\{X^*; \textbf{A}\right\})
&=\underset{\left\| f \right\|_{\mathcal{U}(K)}\le 1}{\text{sup}}\left| \frac{1}{4\pi}\int_Cf(\xi)d\mu(\xi)-\frac{1}{4\pi} \int_Cf_C(\xi)d\omega(\xi)\right|\\
&=\underset{\left\| f \right\|_{\mathcal{U}(K)}\le 1}{\text{sup}}\left| \int_Cf(x)d(\mu-\omega)(\xi)\right|\\
&=\underset{\left\| f \right\|_{\mathcal{U}(K)}\le 1}{\text{sup}}\left| \frac{1}{4\pi}\left\langle f,f_{\mu-\omega} \right\rangle_{\mathcal{U}(K)}\right|=\left\| f_{\mu-\omega} \right\|_{\mathcal{U}(K)}=\left\| \mu-\omega \right\|_{(X,\mathcal{T})}\\
&=D(\omega;\mu).
\end{split}
\end{equation}
Combine ~\eqref{a3d4a22} and ~\eqref{a3d4a2d2}, we have $D_{\min}(\textbf{x};\textbf{A})=\sqrt{\mathcal{E}(\mu(\textbf{x})-\omega(\textbf{x}))}$. This completes the proof.

 This theorem further shows that the maximum value of an integral for potential fields with normalization in a unit can be expressed by the square root of the energy of the charge distribution: $D(0;\mu)=\sqrt{\mathcal{E}(\mu)}$. Some properties of the measures where points with small energy are quantified on how to best approximate the integrals with respect to the equilibrium measure are provided in \cite{damelin2010energy}.

\section{Application}
In a real application, our goal is to generate a large number of points on the sphere that is well-separated and sufficient to cover for the optimal configuration. It is widely used for interpolation and spherical $t$-design \cite{bondarenko2013optimal}. In this section, we will discuss two cases, the first one is to integrate a Franke function for the sphere \cite{renka1988multivariate}, and the second is the distribution of point systems for different kernels to analyze their discrepancies.

\textbf{Integrating a Franke function.}
 The classical kernel approximation of a general function $f:\mathbb{R}^{d+1} \to \mathbb{R}$ can be written as \begin{equation}\label{ad561}
f(\textbf{x})\approx \hat{f}(\textbf{x})=\sum_{i=1}^{N}w_iK(\textbf{x}_i-\textbf{x}'_i).  
\end{equation}
For distinct scattered location $\textbf{x}'_i,i=1,\cdots,N$ with scalar-valued observational data $\textbf{y}= \left[ y_1,\cdots,y_N \right]^T$, a general interpolant can be parameterized as
\begin{equation}\label{ad562}
f(\textbf{x})\approx \hat{f}(\textbf{x})=\sum_{i=1}^{N}w_iK(\textbf{x}_i-\textbf{x}'_i)+\sum_{j=1}^{M}b_jp_j(\textbf{x}), 
\end{equation}
where kernel $K(\cdot)$ acting on the geodesic distance between the center $\textbf{x}'_i$ and the query direction $\textbf{x}$, $p_1(x),\cdots,p_M(x)$ forms a basis for the $M= \mathrm{C}_{s+m'-1}^{m'-1}$-dimensions linear space $ \mathbb{R}_{m-1}^s$ of polynomials of total degree less than or equal to $m'-1$ in $s$ variables. 

The coefficients $\textbf{w}= \left[ w_1,\cdots,w_N \right]^T$ and $\textbf{b}= \left[ b_1,\cdots,b_N \right]^T$ are solutions to the linear equations
\begin{equation}\label{ad563}
(\textbf{K}+\sigma^2\textbf{I})\cdot \textbf{w}+\textbf{p}\cdot\textbf{b}=\textbf{y}.  
\end{equation}
Since enforcing the interpolation condition $f(\textbf{x})\approx \hat{f}(\textbf{x})$, leads to a system of $N$ linear equation with $N+M$ unknown coefficient $w_i$ and $b_j$, and
\begin{equation}\label{ad564}
\sum_{j=1}^{M}w_ip_j(\textbf{x})=0. \ \  j=1,\cdots, M, 
\end{equation}
where $\textbf{K}$ is a matrix with the component $K_{ij}(\cdot)=K( \left\| \textbf{x}_i-\textbf{x}_j \right\|_2)$, $\sigma$ is a smoothing parameter that controls the approximation of the target $f(\textbf{x})$ to fit the observations $\textbf{y}$. 

If $\textbf{K}$ is positive definite and $\textbf{p}$ has a full column rank, the solution for  $\textbf{w}$ and $\textbf{b}$ would be unique. If the chosen $\textbf{K}$  is conditionally positive definite of order $m'$ and $\textbf{p}$ has a full column rank, the solution would be uniquely provided that the degree of the monomial terms is at least $m'-1$ \cite{fasshauer2007meshfree,wahba1990spline}.

Here, our goal is to integrate scattered observations of the Franke function with smoothed parameters for the sphere \cite{renka1988multivariate} defined by
\begin{equation}
\label{nined1s531}
\begin{split}
f(x,y,z):=&\frac{3}{4}\exp(-\frac{(9x-2)^2}{4}-\frac{(9y-2)^2}{4}-\frac{(9z-2)^2}{4})\\
&+\frac{3}{4}\exp(-\frac{(9x+1)^2}{49}-\frac{(9y+1)^2}{10}-\frac{(9z+1)^2}{10})\\
&+\frac{1}{2}\exp(-\frac{(9x-7)^2}{4}-\frac{(9y-3)^2}{4}-\frac{(9z-5)^2}{4})\\
&-\frac{1}{5}\exp(-\frac{(9x-4)^2}{4}-(9y-7)^2-(9z-5)^2),\ \ \ (x,y,z)^T\in\mathbb{S}^2.
\end{split}
\end{equation}
Here, we consider Pycke's statistic $K_A(\eta_i,\eta_j)=-\frac{1}{4 \pi}\ln\frac{e}{2}(1-\eta_i \cdot \eta_j)$ to interpolate the target with the cases of its first $K_A^{(1)}(\eta_i,\eta_j)$ and second $K_A^{(2)}(\eta_i,\eta_j)$ order of derivatives, respectively.

Inspired by \cite{2023submitted}, in order to scale the evaluations with unique points, suppose that we have already generated $n$ points the interpolation points are generated sequentially by
\begin{equation}\label{ad5643s}
\eta_{n+1}=\underset{\eta\in \mathbb{S}^2}{\text{arg}\min}\sum_{i=1}^{n}K(\eta_i,\eta).
\end{equation}
The initial point we choose
\begin{equation}\label{ad56433s}
\eta_{1}=\underset{\eta\in \mathbb{S}^2}{\text{arg} \max} \phi(\eta-\eta_i),
\end{equation}
where $\phi(\cdot)$ follows Gaussian distribution.
Thus, the formulas to generate interpolation points can be written as follows.
\begin{equation}
\label{nined1s531a}
\begin{split}
&\eta_{n+1}=\underset{\eta\in \mathbb{S}^2}{\text{arg}\min}\sum_{i=1}^{n}-\frac{1}{4 \pi}\ln\frac{e}{2}(1-\eta_i \cdot \eta),\\
&\eta_{n+1}^{(1)}=\underset{\eta\in \mathbb{S}^2}{\text{arg}\min}\sum_{i=1}^{n}\frac{1}{1-\eta_i \cdot \eta},\\
&\eta_{n+1}^{(2)}=\underset{\eta\in \mathbb{S}^2}{\text{arg}\min}\sum_{i=1}^{n}\frac{1}{(1-\eta_i \cdot \eta)^2}.
\end{split}
\end{equation}
We use a spherical coordinate system $(r,\theta,\varphi)$, $r=1$ represents the radial distance is equal to $1$ for our experiment, polar angle $\theta\in\left[ 0,\pi \right]$ represents the angle with respect to the polar axis, azimuthal angle $\varphi\in[0,2\pi)$ represents the angle of rotation from the initial meridian plane. The Cartesian coordinates can be retrieved from the spherical coordinate by
\begin{equation}
\label{nined1s5312a}
\begin{split}
&x=\sin\theta \cos\varphi,\\
&y=\sin\theta \sin\varphi,\\
&z= \cos\theta.
\end{split}
\end{equation}
%误差，参数分析
%https://github.com/daixiongming/Operator-free-equidistribution-on-the-sphere/blob/main/scripts/interpolate.b.py
\begin{figure}[!htb]
\minipage{0.32\textwidth}
  \includegraphics[width=\linewidth]{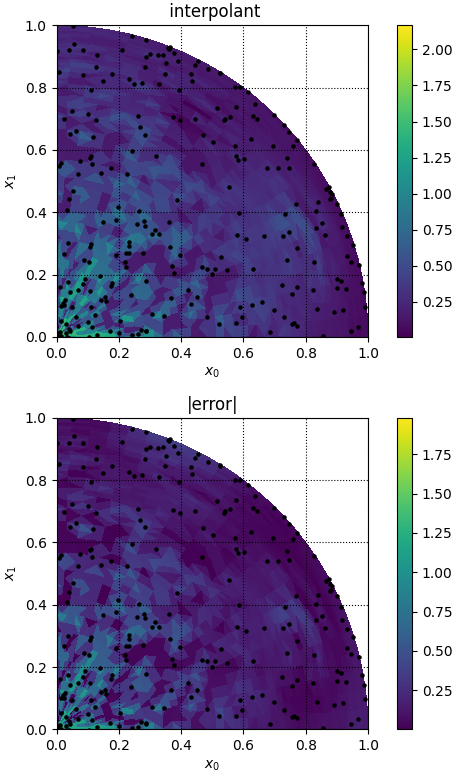}
  %\caption{1.The discrepancy is the average error 0.43376 for 1000 nodes here. }\label{fig:awesome_image11}
\endminipage\hfill
\minipage{0.32\textwidth}
  \includegraphics[width=\linewidth]{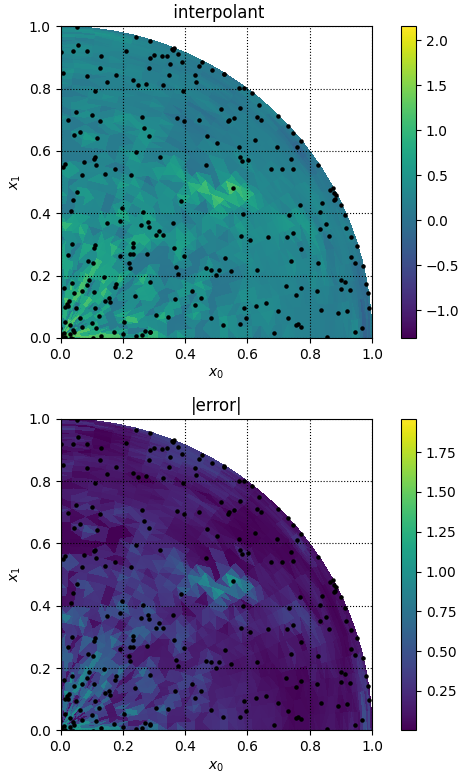}
  %\caption{2.The discrepancy is the average error 0.43659 for 1000 nodes here.}\label{fig:awesome_image21}
\endminipage\hfill
\minipage{0.32\textwidth}%
  \includegraphics[width=\linewidth]{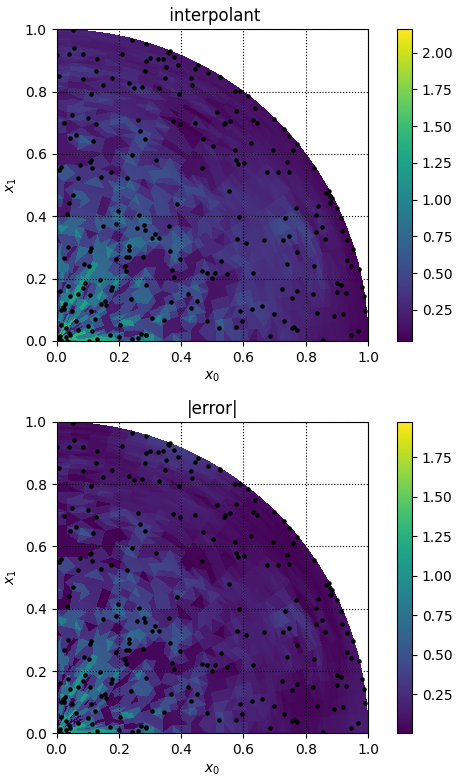}
  %\caption{3.The discrepancy is the average error 0.43361 for 1000 nodes here.}\label{fig:awesome_image31}
\endminipage
\caption{Partition of unity property of the interpolant and the corresponding error distribution. Left: The discrepancy is the average error of 0.43376 for $K_A(\eta_i,\eta_j)$ with 1000 nodes\cite{cui1997equidistribution}. Middle: The discrepancy is the average error of 0.43659 for $K_A^{(1)}(\eta_i,\eta_j)$ with 1000 nodes. Right: The discrepancy is the average error of 0.43361 for $K_A^{(2)}(\eta_i,\eta_j)$ with 1000 nodes.}\label{fig:awesome_image313}
\end{figure}
\begin{table*}[ht]  
  \centering 
%\begin{center}
\caption{The generalized discrepancy of integrated nodes.} 
 \label{tab:data11} 
\begin{tabular}{ | c | c | c |c | }
  \hline
  \# of points & $D(\left\{\eta_1,\cdots,\eta_N\right\};K_A)$ \cite{cui1997equidistribution} & $D(\left\{\eta_1,\cdots,\eta_N\right\};K_A^{(1)})$ & $D(\left\{\eta_1,\cdots,\eta_N\right\};K_A^{(2)})$ \\ \hline
  $15$ & 0.68137655 & 0.68339536 & 0.69559213\\ \hline
  $43$ & 0.59310219 & 0.59549629 & 0.61001457 \\ \hline
  $86$ & 0.54524042 & 0.54779668 & 0.56333339 \\ \hline
  $151$ & 0.51181878 & 0.51446924 & 0.53060611 \\ \hline
  $206$ & 0.49792697 & 0.50061168 & 0.51696945 \\ \hline
    $313$ & 0.47840121 & 0.48112904 & 0.47823923 \\ \hline
  $529$ & 0.45436048 & 0.45713295 & 0.45419589 \\ \hline
    $719$ & 0.44430551 & 0.44709377 & 0.44413999 \\ \hline
  $998$ & 0.43388233 & 0.43668512 & 0.43371597 \\ \hline
\end{tabular}
%\end{center}
\end{table*} 
The Cartesian coordinate of the point on the sphere $\eta=(x,y,z)$. We plot the interpolant in spherical coordinates under three different kernel interpolations in \autoref{fig:awesome_image313}. It shows the point system with the minimum discrepancy is by $K_A^{(2)}(\eta_i,\eta_j)$.

%分析那个误差和discrepancy:
%分析那个表格：
Table~\ref{tab:data11} provides the computed values of the generalized discrepancy for different kernels from the same family. Among them, the best point system is from the second order of derivatives.
%https://github.com/daixiongming/Operator-free-equidistribution-on-the-sphere/blob/main/scripts/interpolate.b.py
\begin{figure}[!htb]
\minipage{0.32\textwidth}
  \includegraphics[width=\linewidth]{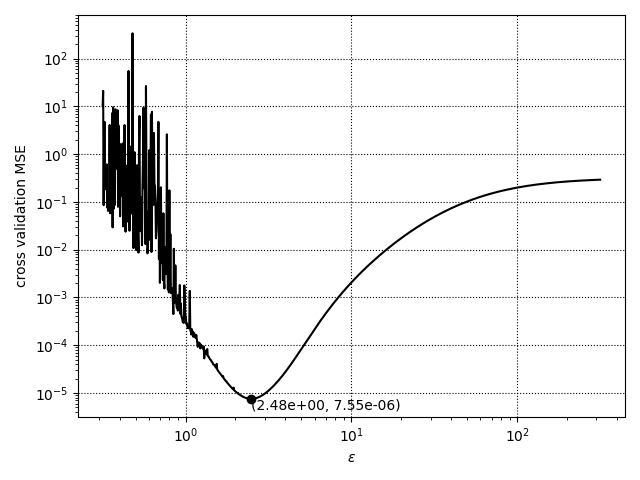}
 % \caption{A really Awesome Image}\label{fig:awesome_image122}
\endminipage\hfill
\minipage{0.32\textwidth}
  \includegraphics[width=\linewidth]{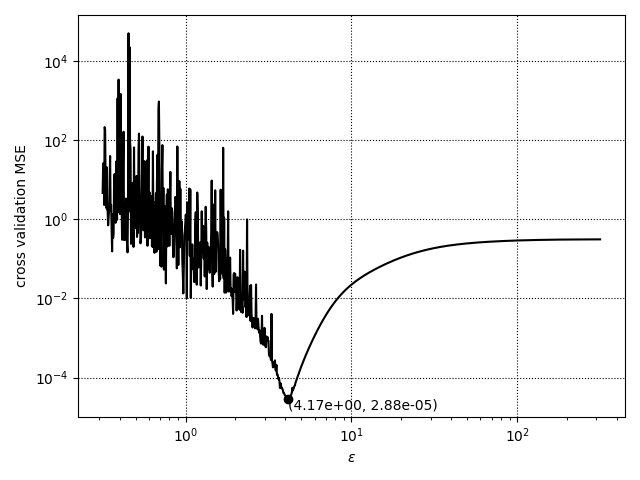}
  %\caption{A really Awesome Image}\label{fig:awesome_image222}
\endminipage\hfill
\minipage{0.32\textwidth}%
  \includegraphics[width=\linewidth]{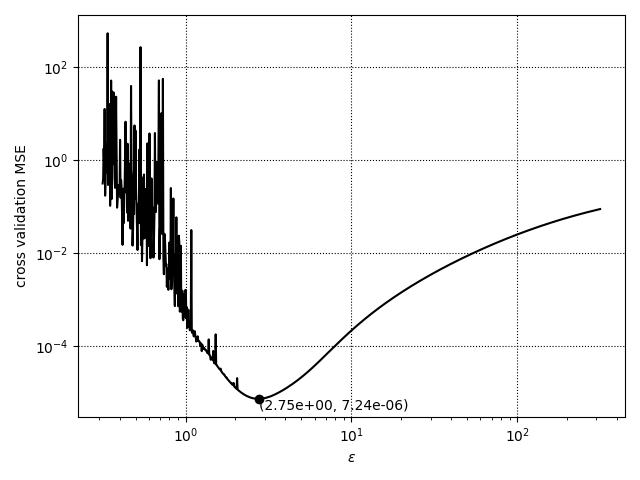}
  %\caption{A really Awesome Image}\label{fig:awesome_image322}
\endminipage
  \caption{Partition of unity property of the interpolant, as a function of the kernel parameter $\varepsilon$ in $d=3$ for $K_A(\eta_i,\eta_j)$\cite{cui1997equidistribution} (left), \ $K_A^{(1)}(\eta_i,\eta_j)$ (middle) and $K_A^{(2)}(\eta_i,\eta_j)$ (right).}\label{fig:awesome_image322}

\end{figure}
%分析那个4到6图的一些特点：

%crossvalidation 分析：
We further estimate the kernel parameter $\varepsilon$ with $N=1000$, by minimizing the mean square error for a fit to the data based on an interpolant.

%Estimating the smoothing parameter: page 164 \cite{wahba1990spline}
%page 55
From ~\eqref{ad562}, ~\eqref{ad563} and ~\eqref{ad564}, the coefficient vector $\textbf{w}= \left[ w_1,\cdots,w_N \right]^T$ and $\textbf{y}= \left[ y_1,\cdots,y_N \right]^T$ are determined by interpolating the observational data $\textbf{y}= \left[ y_1,\cdots,y_N \right]^T$.
\begin{equation}\label{ad56433s1}
f(x_i)=y_i, i=1,\cdots,N,
\end{equation}
which is equivalent to solving the linear system $\textbf{c}=\left[\textbf{w},\textbf{b}\right]^T$,
\begin{equation*}
 Q\textbf{c}=\textbf{y}, \ \ \ Q=g(K(\left\| x_i-x_j \right\|)),
\end{equation*}
where $g(x)$ is a function of $x$.

Inspired by \cite{rippa1999algorithm}, let $U^{(v)}$ the subset obtained by removing the point $x_v$ from $U$ and by $\textbf{y}^{(v)}=\left[ y_1^{(v)},\cdots,y_{v-1}^{(v)},y_{v+1}^{(v)},\cdots,y_N^{(v)} \right]^T$ the vector obtained by removing the element $y_v$ from $\textbf{y}$.

From the perspective of the interpolant
\begin{equation}\label{ad56433s2}
f^{(k)}(x)=\sum_{j=1,j\neq v}^{N}w_j^{(v)}g(K(\left\| x_j-x \right\|)),
\end{equation}
where $\textbf{a}^{(v)}=\left[ w_1^{(v)},\cdots,w_{v-1}^{(v)},w_{v+1}^{(v)},\cdots,w_N^{(v)} \right]^T$ is determined by the interpolation conditions
\begin{equation*}
f^{(k)}(x_i)=y_i, i=1,\cdots,N, i\neq v.
\end{equation*}
which is equivalent to solving
\begin{equation}\label{ad56433s3d}
Q^{(v)}w^{(v)}=f^{(v)},
\end{equation}
where $Q^{(v)}$ is obtained from $Q$ be removing the $v$-th row and $v$-th column, we can obtain the $v$-th error term by
\begin{equation}\label{ad56433s3}
\varepsilon_v=y_v-f^{v}(x_k).
\end{equation}
As the linear system ~\eqref{ad56433s3d} is of order $(N-1)\times (N-1)$, the time complexity is of order $ \mathcal{O}(N^4)$ for the lower-upper decomposition. Fortunately, in real applications, these error components can be simplified to
\begin{equation}\label{ad56433s31}
\varepsilon_v=\frac{c_v}{G_{vv}^{-1}}.
\end{equation}
where $c_v$ is the $v$-th coefficient in the interpolant $f_i$ based on the full dataset, and $G_{vv}^{-1}$ is the $v$-th element of the inverse of the corresponding interpolant matrix, since the complexity of both $c_v$ and $G_{vv}^{-1}$ is $ \mathcal{O}(N^3)$, the computational load will be scaled greatly \cite{rippa1999algorithm}.

In \autoref{fig:awesome_image322}, the optimal $\varepsilon=2.48$ with the minimum mean square of error $7.55*10^{-6}$ is for $K_A(\eta_i,\eta_j)$. The point systems generated from the first order of derivative become worse with the MSE of $2.88*10^{-5}$ when $\varepsilon=4.17$. While the second order of the derivative is the best at $\varepsilon=2.75$ with the MSE of $7.24*10^{-6}$.

\textbf{Point systems for different kernels on the sphere.} Vlasiuk proposes an algorithm to generate high-dimensional points by a combination of quasi-Monte Carlo methods and weighted Riesz energy minimization embedding with a nearest-neighbor distance function \cite{vlasiuk2018fast}. For the node generation on a unit sphere, we simplify the process from the random sampling and normalize it to project on the sphere, which ensures the node must be restricted to a certain compact set $ \mathbb{S}^2$. The schema can be described as follows.

(a) 3D nodes are generated Randomly and normalized to ensure that they are within the unit sphere.

(b) Set up $K'$ nearest neighbors of each node $r=\left\| x-x_i \right\|$. 

(c) Compute the Riesz weight for each node from the corresponding $r$ and normalize it.

(d) Sum the entire weights and find the mean as the discrepancy by $D=\frac{1}{N^2}\sum_{i=1}^{N}\sum_{j=1}^{N}K(\left\| x_i-x_j \right\|)$.

(e) Perform $T$ iterations of the partial gradient descent on the Cui-Freeden discrepancy kernel $K=2-2\log(1+\frac{r}{2})$. Let the configuration by $t$th iteration is $x_i^{t}$, we have $x_i^{0}=x_i,i=1,\cdots,N$, $N$ denotes the number of nodes. Given a node $x_i^{(t)}$ with $K'$ nearest neighbors $x_{j(i,k)}^{(t)}, k=1,\cdots,K'$, the weighted vector sum is
\begin{equation}\label{ad56433s3d1}
g_i^{(t)}=s\sum_{k=1}^{K'}\frac{x_i^{(t)}-x_{j(i,k)}^{(t)}}{\left\| x_i^{(t)}-x_{j(i,k)}^{(t)} \right\|^{s+2}}, 1\le i\le N,
\end{equation}
and the neighbor indices $j(i,k)$ will be updated after a few iterations. The $t+1$th iteration node can be written as
\begin{equation}\label{ad1}
x_i^{(t+1)}=x_i^{(t)}+\frac{\Delta(x_i^{(t)})}{t+C_2}\frac{g_i^{(t)}}{\left\| g_i^{(t)} \right\|}, x_i \in \mathbb{S}^2,
\end{equation}
where $C_2=19$ denotes a fixed offset to control the step size between $x_i^{(t)}$ and $x_i^{(t+1)}$. 
%https://github.com/daixiongming/Operator-free-equidistribution-on-the-sphere/tree/main/BRieszk-master/src
\begin{figure}[!htb]
\minipage{0.32\textwidth}
  \includegraphics[width=\linewidth]{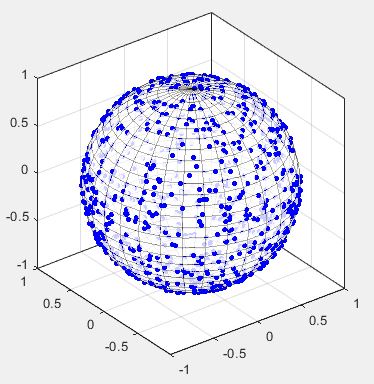}
  %\caption{The discrepancy is the average summation of the kernel 0.092856 for 1000 nodes here.}\label{fig:awesome_image1}
\endminipage\hfill
\minipage{0.32\textwidth}
  \includegraphics[width=\linewidth]{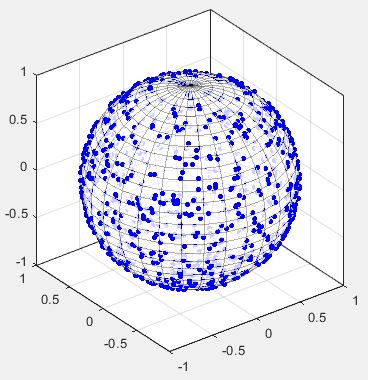}
  %\caption{The discrepancy is the average summation of the kernel 0.051135 for 1000 nodes here.}\label{fig:awesome_image2}
\endminipage\hfill
\minipage{0.32\textwidth}%
  \includegraphics[width=\linewidth]{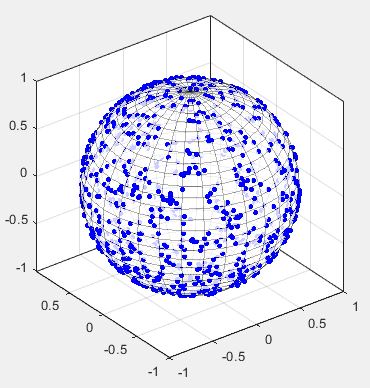}
  %\caption{The discrepancy is the average summation of the kernel 0.015885 for 1000 nodes here.}\label{fig:awesome_image3}
\endminipage
\caption{The node rendering distribution generated from different kernels on the sphere. Left: The discrepancy is 0.092856 for $K$\cite{cui1997equidistribution}. Middle: The discrepancy is 0.051135 for $K^{(1)}$. Right: The discrepancy is 0.015885 for $K^{(2)}$.}\label{fig:awesome_image1-3}
\end{figure}

\begin{table*}[ht]  
  \centering 
%\begin{center}
\caption{The generalized discrepancy of discretized nodes.} 
 \label{tab:data22} 
\begin{tabular}{ | c | c | c |c | }
  \hline
  \# of points & $D(\left\{x_1,\cdots,x_N\right\};K)$\text{\cite{cui1997equidistribution}} & $D(\left\{x_1,\cdots,x_N\right\};K^{(1)})$ & $D(\left\{x_1,\cdots,x_N\right\};K^{(2)})$ \\ \hline
  $15$ & 0.26549542 & 0.19032185 & 0.09912521 \\ \hline
  $43$ & 0.17015372 & 0.10908718 & 0.04564855 \\ \hline
  $86$ & 0.13308495 & 0.08021861 & 0.02974657 \\ \hline
  $151$ & 0.1125596 & 0.06505319 & 0.02221502 \\ \hline
  $206$ & 0.1074984 & 0.06141475 & 0.02050324 \\ \hline
  $313$ & 0.0989759 & 0.05538551 & 0.01775417 \\ \hline
  $529$ & 0.0882991 & 0.04801519 & 0.01455116 \\ \hline
  $719$ & 0.0864768 & 0.04677885 & 0.01403181 \\ \hline
  $998$ & 0.0844954 & 0.04544187 & 0.01347624 \\ \hline
\end{tabular}
%\end{center}
\end{table*} 
%write the Table 2 analysis!

\autoref{fig:awesome_image1-3} shows the node distribution from different kernels, the discrepancy is calculated by the average summation of the kernels for different point distributions on the sphere with 1000 nodes, and the point system generated from the second order of derivatives has the minimum discrepancy.

Table~\ref{tab:data22} provides the computed values of the generalized discrepancy with different numbers of nodes for different kernels from the same family. Among them, the best point system is from the second order of derivatives.

\section{Conclusion}
Generating equidistributed pointsets on the sphere is practical of importance, which generally involves pseudodifferential operators and Beltrami operators to give a quantifying criterion,  it limits to the kernel-self when there exists a closed-form expression. We use the advantage of Legendre's ODE and further explore latent point systems within error bounds, We consider the properties of the kernel with continuity and derivative, Legendre's ODE and spherical harmonic theoretics to find a new criterion of equidistributed pointsets where the discrepancy becomes smaller, and propose a generalized minimum discrepancy. Our kernel-derivative model can explore latent potential point systems that have the minimum discrepancy with operators-free, which has been verified by several quantitive tests in our experiments.

%Step 1: Restate the problem
%Step 2: Sum up the paper
%Step 3: Discuss the implications
%Research paper conclusion examples
%Frequently asked questions about research paper conclusions
\section*{Acknowledgments}
This was supported in part by BRBytes project.

%Bibliography
\bibliographystyle{unsrt}  
\bibliography{references}

\begin{thebibliography}{10}

\bibitem{hlawka1982gleichverteilung}
Edmund Hlawka.
\newblock Gleichverteilung auf produkten von sph{\"a}ren.
\newblock 1982.

\bibitem{grabner1991erdos}
Peter~J Grabner.
\newblock Erd{\"o}s-tur{\'a}n type discrepancy bounds.
\newblock {\em Monatshefte f{\"u}r Mathematik}, 111:127--135, 1991.

\bibitem{grabner1993spherical}
Peter~J Grabner and Robert~F Tichy.
\newblock Spherical designs, discrepancy and numerical integration.
\newblock {\em mathematics of computation}, 60(201):327--336, 1993.

\bibitem{rakhmanov1994minimal}
Evguenii~A Rakhmanov, Edward~B Saff, and YM1306011 Zhou.
\newblock Minimal discrete energy on the sphere.
\newblock {\em Mathematical Research Letters}, 1(6):647--662, 1994.

\bibitem{rakhmanov1995electrons}
Evguenii~A Rakhmanov, EB~Saff, and YM~Zhou.
\newblock Electrons on the sphere.
\newblock In {\em COMPUTATIONAL METHODS AND FUNCTION THEORY 1994}, pages 293--309. World Scientific, 1995.

\bibitem{kuijlaars1998asymptotics}
Arno Kuijlaars and E~Saff.
\newblock Asymptotics for minimal discrete energy on the sphere.
\newblock {\em Transactions of the American Mathematical Society}, 350(2):523--538, 1998.

\bibitem{saff1997distributing}
Edward~B Saff and Amo~BJ Kuijlaars.
\newblock Distributing many points on a sphere.
\newblock {\em The mathematical intelligencer}, 19:5--11, 1997.

\bibitem{andrievskii1999discrepancy}
Vladimir~V Andrievskii, H-P Blatt, and M~Go{\"u}tz.
\newblock Discrepancy estimates on the sphere.
\newblock {\em Monatshefte f{\"u}r Mathematik}, 128:179--188, 1999.

\bibitem{gotz2001note}
Mario G{\"o}tz and Edward~B Saff.
\newblock Note on d—extremal configurations for the sphere in ℝ d+ 1.
\newblock In {\em Recent Progress in Multivariate Approximation: 4th International Conference, Witten-Bommerholz (Germany), September 2000}, pages 159--162. Springer, 2001.

\bibitem{bajnok2015constructive}
B{\'e}la Bajnok, Steven~B Damelin, Jenny Li, and Gary~L Mullen.
\newblock A constructive finite field method for scattering points on the surface of $ d $-dimensional spheres.
\newblock {\em arXiv preprint arXiv:1512.02984}, 2015.

\bibitem{damelin2003energy}
Steven~B Damelin and Peter~J Grabner.
\newblock Energy functionals, numerical integration and asymptotic equidistribution on the sphere.
\newblock {\em Journal of Complexity}, 19(3):231--246, 2003.

\bibitem{sloan2004extremal}
Ian~H Sloan and Robert~S Womersley.
\newblock Extremal systems of points and numerical integration on the sphere.
\newblock {\em Advances in Computational Mathematics}, 21:107--125, 2004.

\bibitem{hesse2006lower}
Kerstin Hesse.
\newblock A lower bound for the worst-case cubature error on spheres of arbitrary dimension.
\newblock {\em Numerische Mathematik}, 103:413--433, 2006.

\bibitem{narcowich2010leveque}
Francis~J Narcowich, Xingping Sun, Joseph~D Ward, and Zongmin Wu.
\newblock Leveque type inequalities and discrepancy estimates for minimal energy configurations on spheres.
\newblock {\em Journal of Approximation Theory}, 162(6):1256--1278, 2010.

\bibitem{brauchart2008optimal}
Johann Brauchart.
\newblock Optimal logarithmic energy points on the unit sphere.
\newblock {\em Mathematics of Computation}, 77(263):1599--1613, 2008.

\bibitem{qiu2000numerical}
Yiqi Qiu, DM~Sloan, and Tao Tang.
\newblock Numerical solution of a singularly perturbed two-point boundary value problem using equidistribution: analysis of convergence.
\newblock {\em Journal of computational and applied mathematics}, 116(1):121--143, 2000.

\bibitem{nigrini2007benford}
Mark~J Nigrini and Steven~J Miller.
\newblock Benford’s law applied to hydrology data—results and relevance to other geophysical data.
\newblock {\em Mathematical Geology}, 39:469--490, 2007.

\bibitem{gustafsson2000reconstructing}
Bj{\"o}rn Gustafsson, Chiyu He, Peyman Milanfar, and Mihai Putinar.
\newblock Reconstructing planar domains from their moments.
\newblock {\em Inverse Problems}, 16(4):1053, 2000.

\bibitem{druskin2010adaptive}
Vladimir Druskin, Chad Lieberman, and Mikhail Zaslavsky.
\newblock On adaptive choice of shifts in rational krylov subspace reduction of evolutionary problems.
\newblock {\em SIAM Journal on Scientific Computing}, 32(5):2485--2496, 2010.

\bibitem{cui1997equidistribution}
Jianjun Cui and Willi Freeden.
\newblock Equidistribution on the sphere.
\newblock {\em SIAM Journal on Scientific Computing}, 18(2):595--609, 1997.

\bibitem{freeden1999constructive}
Willi Freeden and Volker Michel.
\newblock Constructive approximation and numerical methods in geodetic research today--an attempt at a categorization based on an uncertainty principle.
\newblock {\em Journal of Geodesy}, 73:452--465, 1999.

\bibitem{freeden2017integration}
Willi Freeden and Martin Gutting.
\newblock {\em Integration and cubature methods: A geomathematically oriented course}.
\newblock CRC Press, 2017.

\bibitem{watson1967another}
Geoffrey~S Watson.
\newblock Another test for the uniformity of a circular distribution.
\newblock {\em Biometrika}, 54(3-4):675--677, 1967.

\bibitem{beran1968testing}
RJ~Beran.
\newblock Testing for uniformity on a compact homogeneous space.
\newblock {\em Journal of Applied Probability}, 5(1):177--195, 1968.

\bibitem{gine1975invariant}
Evarist Gin{\'e}.
\newblock Invariant tests for uniformity on compact riemannian manifolds based on sobolev norms.
\newblock {\em The Annals of statistics}, 3(6):1243--1266, 1975.

\bibitem{pycke2007decomposition}
Jean-Renaud Pycke.
\newblock A decomposition for invariant tests of uniformity on the sphere.
\newblock {\em Proceedings of the American Mathematical Society}, 135(9):2983--2993, 2007.

\bibitem{pycke2007u}
J-R Pycke.
\newblock U-statistics based on the green's function of the laplacian on the circle and the sphere.
\newblock {\em Statistics \& probability letters}, 77(9):863--872, 2007.

\bibitem{prentice1978invariant}
MJ~Prentice.
\newblock On invariant tests of uniformity for directions and orientations.
\newblock {\em The Annals of Statistics}, pages 169--176, 1978.

\bibitem{freeden1980integral}
Willi Freeden.
\newblock On integral formulas of the (unit) sphere and their application to numerical computation of integrals.
\newblock 1980.

\bibitem{muller2006spherical}
Claus M{\"u}ller.
\newblock {\em Spherical harmonics}, volume~17.
\newblock Springer, 2006.

\bibitem{ruzhansky2009pseudo}
Michael Ruzhansky and Ville Turunen.
\newblock {\em Pseudo-differential operators and symmetries: background analysis and advanced topics}, volume~2.
\newblock Springer Science \& Business Media, 2009.

\bibitem{gronwall1914degree}
Thomas~Hakon Gronwall.
\newblock On the degree of convergence of laplace’s series.
\newblock {\em Transactions of the American Mathematical Society}, 15(1):1--30, 1914.

\bibitem{arfken1972mathematical}
George~B Arfken and Hans-Jurgen Weber.
\newblock Mathematical methods for physicists.
\newblock 1972.

\bibitem{choirat2013computational}
Christine Choirat and Raffaello Seri.
\newblock Computational aspects of cui-freeden statistics for equidistribution on the sphere.
\newblock {\em Mathematics of Computation}, 82(284):2137--2156, 2013.

\bibitem{gradshteyn2014table}
Izrail~Solomonovich Gradshteyn and Iosif~Moiseevich Ryzhik.
\newblock {\em Table of integrals, series, and products}.
\newblock Academic press, 2014.

\bibitem{damelin2008walk}
Steven~B Damelin.
\newblock A walk through energy, discrepancy, numerical integration and group invariant measures on measurable subsets of euclidean space.
\newblock {\em Numerical Algorithms}, 48(1-3):213--235, 2008.

\bibitem{damelin2010energy}
Steven~B Damelin, Fred~J Hickernell, David~L Ragozin, and Xiaoyan Zeng.
\newblock On energy, discrepancy and group invariant measures on measurable subsets of euclidean space.
\newblock {\em Journal of Fourier Analysis and Applications}, 16(6):813--839, 2010.

\bibitem{bondarenko2013optimal}
Andriy Bondarenko, Danylo Radchenko, and Maryna Viazovska.
\newblock Optimal asymptotic bounds for spherical designs.
\newblock {\em Annals of mathematics}, pages 443--452, 2013.

\bibitem{renka1988multivariate}
Robert~J Renka.
\newblock Multivariate interpolation of large sets of scattered data.
\newblock {\em ACM Transactions on Mathematical Software (TOMS)}, 14(2):139--148, 1988.

\bibitem{fasshauer2007meshfree}
Gregory~E Fasshauer.
\newblock {\em Meshfree approximation methods with MATLAB}, volume~6.
\newblock World Scientific, 2007.

\bibitem{wahba1990spline}
Grace Wahba.
\newblock {\em Spline models for observational data}.
\newblock SIAM, 1990.

\bibitem{2023submitted}
Blind Reviewing.
\newblock Weighted riesz particles.
\newblock {\em NIPS2023 submitted}, 360(3):1559--1580, 2023.

\bibitem{rippa1999algorithm}
Shmuel Rippa.
\newblock An algorithm for selecting a good value for the parameter c in radial basis function interpolation.
\newblock {\em Advances in Computational Mathematics}, 11:193--210, 1999.

\bibitem{vlasiuk2018fast}
O~Vlasiuk, Timothy Michaels, Natasha Flyer, and Bengt Fornberg.
\newblock Fast high-dimensional node generation with variable density.
\newblock {\em Computers \& Mathematics with Applications}, 76(7):1739--1757, 2018.

\end{thebibliography}

\end{document}